\documentclass[11pt]{article}   

\usepackage{amssymb,amsfonts,amsmath,bm}
\usepackage{epsfig}
\usepackage{subfigure}
\usepackage{bbm,color}
\usepackage{hyperref}
\usepackage{graphicx}

 \usepackage{mathptmx}     

\usepackage{pstricks}
 \usepackage{pgfplots}

\newcommand{\rme}{{\rm e}}
\newcommand{\rmd}{{\rm d}}

\newcommand{\la}{{\lambda}}
\newcommand{\ep}{{\varepsilon}}
\newcommand{\1}{\mathbbm{1}}
\newcommand{\pd}{{\partial}}

\newcommand{\PP}{{\mathbb P}}

\newcommand{\RR}{{\mathbb R}}

\newcommand{\EE}{{\mathbb E}}
\newcommand{\TT}{{\mathbb A}}

\newcommand{\sL}{{\mathcal L}}
\newcommand{\sP}{{\mathcal P}}
\newcommand{\sF}{{\mathcal F}}

\newtheorem{rem}{Remark}[section]

\newtheorem{ex}{Example}[section]

\newcounter{lis}

\begin{document}
\date{}

\title{First-passage probabilities and invariant distributions of Kac-Ornstein-Uhlenbeck processes}
\author{Nikita Ratanov\\
 Chelyabinsk State University, Russia\\
nikita.ratanov@csu.ru,\\
}
\maketitle

\begin{abstract}
In this talk, we present Ornstein-Uhlenbeck processes 
whose parameters are modulated by an external underlying two-state Markov process. 
Conditional mean  of the such process 
for a given modulation follows the analogue of the Langevin equation controlled by a pair of telegraph processes.
These processes are investigated from the point of view of the first passage probabilities and invariant measures. 

It also examines the {limiting behaviour} under scaling conditions similar to Kac's scaling. 
The limiting processes become different classes of ordinary Ornstein-Uhlenbeck processes.
\end{abstract}

\section{Introduction}
The Ornstein-Uhlenbeck process $X$ can be defined as the solution of the stochastic equation
\begin{equation}
\label{eq:OU}
\rmd X^{\mathrm{OU}}(t)=(a-\gamma X^{\mathrm{OU}}(t))\rmd t+b\rmd W(t),\quad t>0,
\end{equation}
with the initial condition
$X^{\mathrm{OU}}(0)=x,$
where $a,\;b$ and $\gamma$ are constants, $b\ge0;$
$W=W(t)$ is a standard Wiener process.

It is widely accepted that Ornstein-Uhlenbeck processes
can be used as an alternative model of Brownian motion 
that better matches the physical data than the pure Wiener process.

The solution of \eqref{eq:OU} can be written explicitly, 
\begin{equation}
\label{sol:OU}
X^{\mathrm{OU}}(t)=\frac{a}{\gamma}
+\left(x-\frac{a}{\gamma}\right)\rme^{-\gamma t}+b\int_0^t\rme^{-\gamma(t-s)}\rmd W(s).
\end{equation}
Basic properties of $X^{\mathrm{OU}}$ can be derived from \eqref{eq:OU}-\eqref{sol:OU}.
In particular, $X(t)$ has a Gaussian distribution with mean
\[
\EE[X(t)]=\frac{a}{\gamma}+\left(x-\frac{a}{\gamma}\right)\rme^{-\gamma t},
\]
and covatiance
\[
\mathrm{Cov}(X(s),\;X(t))=\frac{b^2}{2\gamma}\rme^{-\gamma s}(\rme^{\gamma t}-\rme^{-\gamma t}),\;s\ge t\ge0.
\]
Having originated from physics,
this model is exploited in various application fields 
as an alternative to Wiener process with an average tendency to return, see \cite{Langevin,Maller}.
For instance, the Va\v{s}\'i\v{c}ek interest rate model, \cite{vasicek}, gave rise to widespread
 financial application of this process. 
The same processes are also intensively used for neuronal modelling, see e.g. \cite{pirozzi,ricciardi,BioSys}.

Processes $X^{\mathrm{OU}}$ are studied  
since the seminal paper by Uhlenbeck and Ornstein, \cite{OU}.
Subsequently, similar processes were constructed on the basis of 
a fractional Wiener process, \cite{biagini,kleptsina}, or,
 in general, by a L\'evy process, \cite{B-NS}, and a fractional L\'evy process, \cite{claudia}.
See also \cite{jacobsen96} for rationales and nonstandard interpretations,
and \cite{kleptsina,kutoyants} for statistics.

Recently, some results have appeared on the Ornstein-Uhlenbeck processes with Markov modulation, see
\cite{huang2014,MCAP2020,BioSys,ICSM,xing}. This approach assumes that
equation \eqref{eq:OU} is modified by Markov oscillation of all parameters.  
The motivation to study the model in this context arises from the following observations. 
First, by allowing alternation of the coefficients of the Langevin equation \eqref{eq:OU}, we open up
new good opportunities for applications. 
The second idea is based on the fact that
 the Brownian motion has some properties that are contrary to physical intuition,
  such as infinite total variation of paths and infinite propagation speeds. 
Since a source of stochasticity such as Markov modulation does not possess these disadvantages,
this modification of the model could be fruitful.  

 In general, the most radical approach to overcoming these difficulties
is to
replace the Wiener process with so-called Kac's telegraph process $\mathbb T(t)$, 
\cite{Kac}.
The telegraph process $\mathbb T(t),\;t\ge0,$ describes a steady state motion of a particle 
with alternating at random times velocities, 
\begin{equation}
\label{def:telegraph}
\mathbb T(t)=\int_0^tc_{\ep(s)}\rmd s,
\end{equation}
where $\ep=\ep(s)\in\{0, 1\}$ is a two-state Markov process 
with switching intensity $\la,$ $\la>0.$ 
This mathematical construct is useful on its own,
for example, for 
physically oriented applications, such as the description of
photon transport in a highly scattering medium, or of neutron transport in a reactor, 
\cite{masoliver92,weiss94,weiss02}.
This model has recently been applied to studies of cosmic microwave background
radiation studies based on the hyperbolic heat equation, \cite{leonenko}.
Application of telegraph processes in modelling of financial markets 
and related mathematical innovations are presented in \cite{Quant,KR}.

In this paper, we study
Markov-modulated Ornstein-Uhlenbeck processes where the Wiener process
is replaced by Kac's telegraph process. We call the result of such a replacement
Kac-Ornstein-Uhlen\-beck processes.

This paper concerns the following principal topics. 
First, we are interested in first passage probabilities of the Kac-Ornstein-Uhlenbeck process,
Sections \ref{sec2} and \ref{sec3}.
This subject is related to applications of persistent random walks which are still of interest, 
\cite{hernan,masoliver92,mendez,BioSys,weiss81}. 
Second, we study invariant measures for Markov process 
which is formed by the Kac-Ornstein-Uhlenbeck process $X$ and the underlying state process $\ep,$
$\langle X(t), \ep(t)\rangle,$   Sections \ref{sec2} and \ref{sec5}.
Third, the limit behaviour of $X$ under some traditional parameter scalings is also analysed, see
Section \ref{sec6}.

\section{Model and main objectives }\label{sec2}
On a given probability space $(\Omega, \sF, \PP),$
consider the standard Winer process $W(t),$ $t\ge0,$
 an irreducible continuous-time Markov chain $\ep=\ep(t),\;t\ge0,$ 
 with a finite state space $\{0,\;1,\;\ldots,\; d\}$
and a random variable $x,$ independent of each other.

We define Markov-modulated Ornstein-Uhlenbeck process $M=M(t),$ 
assuming that the parameters $a,\;b$ and $\gamma$ of the Ornstein-Uhlenbeck process
undergo synchronous switching driven by the underlying process $\ep$.

Precisely, let $a_i, b_i$ and $\gamma_i,$ $i\in\{0,\;1,\;\ldots,\; d\},$ be arbitrary constants, $b_i\ge0.$
The process $M=M(t)$ follows the stochastic equation 
\begin{equation*}
\rmd M(t)=\left(a_{\ep(t)}-\gamma_{\ep(t)}M(t)\right)\rmd t+b_{\ep(t)}\rmd W(t),\quad t>0,
\end{equation*}
with the initial condition $M(0)=x.$ This initial value problem is equivalent to the integral equation  
\begin{equation}
\label{eq:KacLangevind}
M(t)=x+\int_0^t\left(a_{\ep(s)}-\gamma_{\ep(s)}M(s)\right)\rmd s+\int_0^tb_{\ep(s)}\rmd W(s),\qquad t\ge0.
\end{equation}

The study of Markov-modulated Ornstein-Uhlenbeck processes has recently begun,
first in \cite{huang2014,xing}, dealing with the transient behaviour of moments 
and some specific scaling of parameters,
and then in \cite{MCAP2020,BioSys,ICSM} in terms of first passage distributions 
and neural modelling applications.

The solution to the integral equation \eqref{eq:KacLangevind} 
can be expressed by means of piecewise deterministic processes.
Let $\Gamma(t)=\int_0^t\gamma_{\ep(s)}\rmd s$ and $\TT(t)=\int_0^ta_{\ep(s)}\rmd s$ 
be two piecewise linear processes based on the common underlying Markov process $\ep.$ 
We define also a process $\widetilde\Gamma(s,t)$ with inverse time by setting 
$\widetilde\Gamma(s,t)=\Gamma(t)-\Gamma(s)=\int_{s}^{t}\gamma_{\ep(u)}\rmd u,\;$ $0\le s\le t.$
The unique solution to \eqref{eq:KacLangevind} is given by
\begin{equation}
\label{eq:solKacLangevind}
M(t)=x\rme^{-\Gamma(t)}+\int_0^t\rme^{-\widetilde\Gamma(s,t)}\rmd \TT(s)
+\int_0^tb_{\ep(s)}\rme^{-\widetilde\Gamma(s,t)}\rmd W(s).
\end{equation}
Conditionally (for given $\{\ep(s)\}_{s\in[0, t]})$ the random variable $M(t)$ is Gaussian with (random) mean
$X(t)=\EE\left(M(t)~|~\{\ep(s)\}_{s\in[0, t]}\right)$ 
and (random) variance \newline$V(t)=\mathrm{Var}\left(M(t)~|~\{\ep(s)\}_{s\in[0, t]}\right),$
\begin{align}
\label{eq:X}
X(t)=&x\rme^{-\Gamma(t)}+\int_0^ta_{\ep(s)}\rme^{-\widetilde\Gamma(s,t)}\rmd s,\\
\label{eq:V}
V(t)=&\int_0^tb_{\ep(s)}^2\rme^{-2\widetilde\Gamma(s, t)}\rmd s.
\end{align}
See \cite[Theorem 2.1]{huang2014}.

In what follows, we assume that the underlying process
$\ep=\ep(t)\in\{0, \;1\},t\ge0,$ is the two-state continuous-time Markov chain 
with transition intensities $\la_0$ and $\la_1,\;\la_0, \la_1>0.$ Let  
\[
\Lambda=\begin{pmatrix}
   -\la_0   &    \la_0\\
   \la_1   &  -\la_1
\end{pmatrix}
\]
and $\Pi=\Pi(t)=\left(\PP\{\ep(t)=j~|~\ep(0)=i\}\right)_{i, j\in\{0, 1\}},$ $t\ge0,$ 
 be the matrix of transition probabilities. It is known that,  see e. g. \cite{ross},
\begin{equation*}\begin{aligned}
\Pi(t)=&\exp(t\Lambda)\\
=&(2\la)^{-1}\begin{pmatrix}
    \la_1+\la_0\rme^{-2\la t}  & \la_0\left(1-\rme^{-2\la t}\right)   \\
     \la_1\left(1-\rme^{-2\la t}\right)  &     \la_0+\la_1\rme^{-2\la t}
\end{pmatrix},
\end{aligned}\end{equation*}
where $2\la=\la_0+\la_1.$
For arbitrary distribution $\vec\pi$ of the initial state $\ep(0),$ the distribution of $\ep(t)$ is given by
$\vec\pi\Pi(t)$ and the limit 
$\vec\pi_*=\vec\pi_*(\la_0,  \la_1)=\lim\limits_{t\to\infty}\vec\pi\Pi(t)$ is given by
\begin{equation}
\label{eq:pi*}
\pi_*(\la_0,\;\la_1)=(2\la)^{-1}(\la_1,\;\la_0).
\end{equation}

The conditional mean $X=X(t)$ of the Markov-modulated process $M$
obeys the integral equation
\begin{equation}
\label{eq:intKacLangevin}
X(t)=x+\int_0^t\left(a_{\ep(s)}-\gamma_{\ep(s)}X(s)\right)\rmd s,\qquad t\ge0.
\end{equation}
The process $X=X(t),\;t\ge0,$ \eqref{eq:intKacLangevin},
 can be viewed as a (non-Gaussian) Ornstein-Uhlenbeck process, controlled 
by two Kac's telegraph processes, $\Gamma(t)$ and $\TT(t)=\int_0^ta_{\ep(s)}\rmd s,\; t\ge0,$
instead of Brownian motion.
We call it the Kac-Ornstein-Uhlenbeck process.

 By virtue of \eqref{eq:X},
the  process $X=X(t)$ 
sequentially follows the two deterministic patterns, $\phi_0$ and $\phi_1,$ 
switching from one to another randomly after exponentially distributed holding times. 
These patterns are 
defined by 
\begin{equation}\label{def:phi}
    \phi_0(t, x) =\rme^{-\gamma_0 t}\left(x+a_0\int_0^t\rme^{\gamma_0s}\rmd s\right),
\qquad
    \phi_1(t, x) =\rme^{-\gamma_1 t}\left(x+a_1\int_0^t\rme^{\gamma_1s}\rmd s\right).
       \end{equation}   
Let $\gamma_0,\;\gamma_1\neq0.$    Then 
 the patterns $\phi_0$ and $\phi_1$ are determined by the function $\phi(t, x)=\rho+(x-\rho)\rme^{-\gamma t}$ 
with two pairs of parameters, $\langle \rho_0, \gamma_0\rangle$ and $\langle \rho_1, \gamma_1\rangle,$
alternating at random times when the underlying process $\ep$ is switched.
Here $\rho_0=a_0/\gamma_0,$ $\rho_1=a_1/\gamma_1.$ If $\gamma=0,$ then $\phi=x+at.$

If both $\gamma$ are positive, we call the model attracting. 
We also consider the attractive-repulsive model,  $\gamma_0\cdot\gamma_1<0.$ 
The case when one of attraction rates is equal to zero, say $\gamma_0>0,\;\gamma_1=0,$
we call non-strict attractive one.

In all cases, $t\to\phi_0(t, x)$ and $t\to\phi_1(t, x),$ satisfy semigroup property.

If $a_0/\gamma_0=a_1/\gamma_1=:\rho,$ then the solution of \eqref{eq:intKacLangevin}
comes down to the exponential telegraph process.  Namely, one can see that in this case,
formula \eqref{eq:X} is simplified as
\begin{equation}\label{eq:X=}
X(t)=\rme^{-\Gamma(t)}\left(x+\rho\int_0^t\gamma_{\ep(s)}\rme^{\Gamma(s)}\rmd s\right)
=\rho+(x-\rho)\exp(-\Gamma(t)).
\end{equation}
The distribution of such a process is well studied, see \cite{KR,JAP51,JAP56,Quant}.
Note that in this case, process $X,$ \eqref{eq:X=}, is time-homogeneous in the sense of \cite[(2.13)]{JAP56}
with rectifying diffeomorphism $\Phi(x)=\log|x-\rho|.$

Assume that $a_0/\gamma_0\ne a_1/\gamma_1.$ 
To be specific, let $\rho_0=a_0/\gamma_0<a_1/\gamma_1=\rho_1.$ 
We are interested in the probabilities of the first passage of a fixed level $y$
by process $X$.

In the case when the parameters $\rho,\;\gamma$ and the variables $x,\; y$ 
satisfy the conditions
\begin{equation}\label{eq:txy>0}
\gamma>0,\quad \frac{x-\rho}{y-\rho}>1\quad\text{or}\quad 
\gamma<0,\quad 0<\frac{x-\rho}{y-\rho}<1,
\end{equation}
the inverse function $t(x, y)=\phi(\cdot, x)^{-1}(y)$ is positive,
\begin{equation}\label{eq:txy}
t(x, y)=\phi(\cdot, x)^{-1}(y)=\frac{1}{\gamma}\log\frac{x-\rho}{y-\rho}>0.
\end{equation}
If for the pair $\langle \rho,\;\gamma\rangle$
condition \eqref{eq:txy>0} is not satisfied, we set $t(x, y)=+\infty$.
By $t_0(x, y)$ and $t_1(x, y)$ we denote the inverse functions to $\phi_0(\cdot, x)$ and $\phi_1(\cdot, x),$ respectively,
which are defined above, \eqref{eq:txy}. The values $t_0(x, y)$ and  $t_1(x, y)$ 
coincide with the shortest time for process $X$ to reach level $y,$ starting from $x,$ and
without switching states.

\subsection{First passage time}
Let $T(x, y)$ be the time when the process $X=X(t)$ first passes through $y,$ starting from $x=X(0),$
\begin{equation}
\label{def:Txy}
T(x, y)=\inf\{t>0:~X(t)=y~|~X(0)=x\},\quad x\ne y.
\end{equation} 
Random variable $T(x, y)$ has an atomic value at $t(x, y),$ $t(x, y)>0,$ if the particle starting at $x$
reaches $y$ without switching.
Furthermore, the distribution of $T(x, y)$ is updated after each state switch.

The following 
distribution identities hold:
\begin{equation}
\label{eq:distributionT}
\begin{aligned}
  \big[ T(x, y)~|~\ep(0)=0\big]  \stackrel{D}{=}t_0(x, y)&\1_{\{\tau^{(0)}>t_0(x, y)\}}\\
  +\big(\tau^{(0)}&+  \left[ T(\phi(\tau^{(0)}, x), y)~|~\ep(0)=1\right]\big)\1_{\{\tau^{(0)}<t_0(x, y)\}}, \\
  \left[ T(x, y)~|~\ep(0)=1\right]  \stackrel{D}{=}t_1(x, y)&\1_{\{\tau^{(1)}>t_1(x, y)\}}\\
 +\big(\tau^{(1)}&+  \left[ T(\phi(\tau^{(1)}, x), y)~|~\ep(0)=0\right]\big)\1_{\{\tau^{(1)}<t_1(x, y)\}}. 
\end{aligned}
\end{equation}
Here, the exponentially distributed random variables $\tau^{(0)}$ and $\tau^{(1)},$ 
$\tau^{(0)}\sim\mathrm{Exp}(\la_0),$ $\tau^{(1)}\sim\mathrm{Exp}(\la_1),$  do not depend on further dynamics;
$[T ~|~\ep(0)=i]$ denotes the conditional distribution of $T$ under the given initial state $\ep(0)=i.$
The first terms on the right-hand sides of equations \eqref{eq:distributionT} are set to zero
if the corresponding $t_i(x, y)$ becomes equal to $+\infty.$

Our first goal is explicit formulae for the Laplace transforms
 \[\ell_0(q, x, y):=\EE_0[\exp(-qT(x, y))],\]\[\ell_1(q, x, y):=\EE_1[\exp(-qT(x, y))].\] 
 Notice that the functions $\ell_0(q, x, y)$ and $\ell_1(q, x, y)$
serve as cumulative distribution functions for the running minimum 
$\underline X_{\rme_q}:=\min\limits_{0\le t\le\rme_q}X(t)$ if $x>y,$
and complementary cumulative distribution function for the running
maximum $\overline X_{\rme_q}:=\max\limits_{0\le t\le\rme_q}X(t)$ if $x<y$.
Indeed, integrating by parts one can see, $i\in\{0, 1\},$
\[
\begin{aligned}
    \ell_i(q, x, y)&= \int\limits_0^\infty\rme^{-qt}\rmd\PP\{T(x, y)<t~|~\ep(0)=i\} 
    =\int\limits_0^\infty q\rme^{-qt}\PP\{T(x, y)<t~|~\ep(0)=i\}\rmd t \\
    &=\PP\{T(x, y)<\rme_q~|~\ep(0)=i\}=\begin{cases}
    \PP\{\underline X_{\rme_q}<y~|~\ep(0)=i\}  & x>y, \\
      \PP\{\overline X_{\rme_q}>y~|~\ep(0)=i\}& x<y,
\end{cases}
\end{aligned}
\]

Due to \eqref{eq:distributionT}, functions $\ell_0$ and $\ell_1$
obey the coupled integral equations,
\begin{equation}
\label{eq:hatQ}
\left\{
\begin{aligned}
\ell_0(q, x, y)&=\rme^{-(q+\la_0)t_0(x, y)} 
+\int_{0}^{t_0(x, y)}\la_0\rme^{-(q+\la_0)\tau}\ell_1(q, \phi_0(\tau, x), y)\rmd\tau,\\
\ell_1(q, x, y)&=\rme^{-(q+\la_1)t_1(x, y)}
+\int_{0}^{ t_1(x, y)}\la_1\rme^{-(q+\la_1)\tau}\ell_0(q, \phi_1(\tau, x), y)\rmd\tau.
\end{aligned}
\right.
\end{equation}
If condition \eqref{eq:txy>0} is not satisfied for a set of parameters $x, y,\;\rho_i, \gamma_i$ , then
$t_i(x, y)=+\infty;$  the first term on the right-hand side of the corresponding equation of \eqref{eq:hatQ} vanishes,
and the next integral is taken over the entire half-line $[0,\;+\infty).$

Differentiating \eqref{eq:hatQ} with respect to $x,$ and then
integrating by parts we get the coupled differential equations:
\begin{equation}
\label{eq:hatQ++diff}
\left\{
\begin{aligned}
(x-\rho_0)\frac{\pd\ell_0}{\pd x}(q, x, y)=&-\beta_0(q)\ell_0(q, x, y)+\beta_0(0)\ell_1(q, x, y),\\
(x-\rho_1)\frac{\pd\ell_1}{\pd x}(q, x, y)=&\beta_1(0)\ell_0(q, x, y)-\beta_1(q)\ell_1(q, x, y).\\   
     \end{aligned}
\right.
\end{equation}
Here  we used the identities
\begin{equation*}
(x-\rho_0)\frac{\pd\phi_0}{\pd x}(\tau, x)\equiv-\frac{1}{\gamma_0}\frac{\pd\phi_0}{\pd \tau}(\tau, x),\qquad
(x-\rho_1)\frac{\pd\phi_1}{\pd x}(\tau, x)\equiv-\frac{1}{\gamma_1}\frac{\pd\phi_1}{\pd \tau}(\tau, x).
\end{equation*}

In the case of non-strict attraction/repulsion, equations 
\eqref{eq:hatQ} and \eqref{eq:hatQ++diff} can be written similarly.

For example, let $\gamma_1=0, a_1>0,\;\gamma_0>0$ and $x<y.$
Hence, the second equation of \eqref{eq:hatQ} hold with $t_1(x, y)=(y-x)/a_1$ and
$\phi_1(\tau, x)=x+a_1\tau.$ Equivalently, the first equation of system \eqref{eq:hatQ++diff} don't change,
while the second equation turns into 
\begin{equation}
\label{eq:2}
a_1\frac{\pd\ell_1}{\pd x}(q, x, y)=-\la_1\ell_0(q, x, y)+(q+\la_1)\ell_1(q, x, y).
\end{equation}

For various combinations of parameters, system \eqref{eq:hatQ++diff} should be considered in different 
domains with different boundary conditions.
In Section \ref{sec3}, we present explicit formulae for $\ell_0(q, x, y)$ and $\ell_1(q, x, y)$ 
with different parameters.

\subsection{Invariant measures}
Our second goal is to study invariant measures for $X.$

Notice that $\Xi(t)=\langle X(t),\ep(t) \rangle\in\mathbb R\times\{0,1\},$ $t\ge0,$ is the Markov process.
Let $\sP(t, \rmd y~|~x)$ be the transition function,
\[
\sP(t, \rmd y~|~x,\;i)=
\left(p_{ij}(t, x; \rmd y)\right)_{i, j\in\{0, 1\}},
\]
where $p_{ij}(t, x; \rmd y)= \PP\{X(t)\in\rmd y,\;\ep(t)=j~|~X(0)=x,\;\ep(0)=i\}.$
Let $P_t$ be the corresponding Markov semigroup, 
 $\vec f\to P_t\vec f,$ 
where  
\[\begin{aligned}
(P_t\vec f)_i(x)&=\EE\left(\vec f(\Xi(t))~\Big|~\Xi(0)=\langle x,\;i\rangle\right)\\
&=\int_{-\infty}^\infty\sP(t, \rmd y~|~x,\;i)\vec f(y),
  \end{aligned}\]
 for any test-function $\vec f=(f_0,\;f_1).$
The infinitesimal generator for the semigroup $P_t$ is determined by
\begin{equation*}
\sL=\begin{pmatrix}
    -\la_0+(a_0-\gamma_0x)\frac{\rmd}{\rmd x}  &  \la_0  \\ \\
    \la_1  &  -\la_1+(a_1-\gamma_1x)\frac{\rmd}{\rmd x}  
\end{pmatrix}.
\end{equation*}
Indeed, let $f_0,\;f_1$ be a pair of test function. By definition, we get
\[\begin{aligned}
&\frac{\EE_0[f_{\ep(t)}(X^x(t))]-f_0(x)}{t}
=\frac{(1-\la_0t)f_0(\phi_0(t, x))+\la_0tf_1(x)-f_0(x)}{t}+o(t)
\\
&=\la_0f_1(x)-\la_0f_0(x)+f_0'(x)\cdot(a_0-\gamma_0x)+o(t), \qquad t\to0,
  \end{aligned}\]
which give the first row of the matrix $\sL$. The second row is obtained similarly. 

We study invariant measures, which are defined as fixed points of the adjoint semigroup $P_t^*.$ 
Let $K,\;K\subset\RR,$ be an invariant set with respect to the time evolution 
$X(t),\;t\ge0,$ \eqref{eq:intKacLangevin}.
The invariant measure $\vec\mu=(\mu_0,\;\mu_1),$ 
 supported on a set $K\subset\RR,$  is defined by the equation
\[
\vec\mu(\rmd y)=\int_{K}\sP(t, \rmd y~|~x)\vec\mu(\rmd x),\qquad y\in K.
\]

When the invariant measure $\vec\mu$ is determined by the probability density function 
$\vec\pi=\vec\pi(x)=(\pi_0(x),\;\pi_1(x)),$ 
this is equivalent to the boundary value problem for the ordinary differential equation, see e.g. \cite{pavliotis},
\begin{equation}
\label{eq:L*pi}
\sL^*\vec\pi(x)=0,\qquad x\in K.
\end{equation}
Here $\sL^*$ is the adjoint operator to the generator $\sL,$ 
and the following assumptions hold:
\[
\pi_0(x)\ge0,\qquad \pi_1(x)\ge0,\quad \forall x\in K
\]
and
\[
\int_{K}\left(\pi_0(x)+\pi_1(x)\right)\rmd x=1.
\]
The existence of the invariant distribution for the process $X$ and its shape 
depends on signs of the parameters $\gamma_0, \gamma_1,$ 
which determine the boundary conditions to equation \eqref{eq:L*pi}.

The explicit form of the adjoint operator $\sL^*$ and the boundary conditions
can be obtained by integrating by parts in
 \[\begin{aligned}
 &\left(\sL\vec f(x), \vec\pi(x)\right)=\\
 \int_{K}&
 \Big(-\la_0f_0(x)+(a_0-\gamma_0x)f_0'(x)+\la_0f_1(x)\Big)\pi_0(x)\rmd x\\
+&\int_{K}
\Big(\la_1f_0(x)-\la_1f_1(x)+(a_1-\gamma_1x)f_1'(x)\Big)\pi_1(x)\rmd x
\end{aligned} \]
 for any test function $\vec f=(f_0(x),\;f_1(x)).$ 
 We have
\begin{equation}\label{eq:Lfpi}
\begin{aligned}
\left(\sL\vec f(x), \vec\pi(x)\right)&\\
=\big[(a_0-\gamma_0x)f_0(x)\pi_0(x)+(a_1-&\gamma_1x)f_1(x)\pi_1(x)\big]|_{x\in\pd K}\\
+\int_{-\infty}^\infty f_0(x)\big[-\la_0\pi_0(x)+&\la_1\pi_1(x)
+\gamma_0\pi_0(x)+(\gamma_0x-a_0)\pi_0'(x)\big]\rmd x\\
+\int_{-\infty}^\infty f_1(x)\big[\la_0\pi_0(x)-&\la_1\pi_1(x)
+\gamma_1\pi_1(x)+(\gamma_1x-a_1)\pi_1'(x)
\big]\rmd x.
\end{aligned}\end{equation}
Therefore, the adjoint operator $\sL^*$ is defined by the matrix
\begin{equation*}
\sL^*=
\begin{pmatrix}
     (\gamma_0-\la_0)+(\gamma_0x-a_0)\frac{\rmd}{\rmd x}&\la_1 \\ \\
      \la_0&(\gamma_1-\la_1)+(\gamma_1x-a_1)\frac{\rmd}{\rmd x}
\end{pmatrix},
\end{equation*}
and the boundary conditions for \eqref{eq:L*pi} are supplied by setting the \textit{non-integral terms}
of \eqref{eq:Lfpi} to be zero.

Precisely, we have the following system:
\begin{equation}
\label{eq:L*system}
\left\{
\begin{aligned}
    (\gamma_0x-a_0)\frac{\rmd\pi_0(x)}{\rmd x}&=\left(\la_0-\gamma_0\right) \pi_0(x)-\la_1\pi_1(x),  \\
   (\gamma_1x-a_1)\frac{\rmd\pi_1(x)}{\rmd x}&=-\la_0\pi_0(x)+\left(\la_1-\gamma_1\right) \pi_1(x), 
\end{aligned}
\right.\qquad x\in K,
\end{equation}
with the boundary conditions
\begin{equation}
\label{eq:L*bc}
(a_0-\gamma_0x)\pi_0(x)|_{\pd K}=0,\qquad(a_1-\gamma_1x)\pi_1(x)|_{\pd K}=0.
\end{equation}

Below, Section \ref{sec5}, we study the invariant measures under different combinations of parameters of the model.
We distinguish two main cases: the attrac\-ting-only dynamics when both $\gamma$ are positive 
and the mixed attrac\-tion-repulsion case. 

The first passage probabilities are explored in the next section. 
\section{First passage probabilities for the Kac-Ornstein-Uhlenbeck process} \label{sec3}

In this section, 
we obtain some explicit formulae for the Laplace transforms $\ell_0$ and $\ell_1$ 
of the first passage time $T(x, y),$ \eqref{def:Txy}. We will consider two different models,
when the paths of $X$ are alternately attracted to the points $\rho_0=a_0/\gamma_0$ 
and $\rho_1=a_1/\gamma_1$, $\gamma_0,\;\gamma_1>0,$
and the case when one level attracts and the other repels, $\gamma_0\cdot\gamma_1<0$.

We will use the following notations:
\[
\xi_0(x)=\frac{x-\rho_0}{\rho_1-\rho_0},\qquad 
\xi_1(x)=1-\xi_0(x)=\frac{\rho_1-x}{\rho_1-\rho_0}
\]
and
\begin{equation}\label{def:b0b1}
b_{0, 1}=\frac12\left(\beta_0+\beta_1\pm\sqrt{(\beta_0-\beta_1)^2+4\la_0\la_1/(\gamma_0\gamma_1)}\right),
\end{equation}
where \[\beta_0=\beta_0(q)=(q+\la_0)/\gamma_0,\;\beta_1=\beta_1(q)=(q+\la_1)/\gamma_1;\]
by $F(b_0, b_1; b_2; \cdot)$ we denote the Gaussian hypergeometric function,
defined by the series 
\begin{equation}
\label{def:F}
F(b_0, b_1; b_2; z)=1+\sum_{n=1}^\infty\frac{(b_0)_n(b_1)_n}{(b_2)_nn!}z^n,
\end{equation}
if one of the following conditions holds:
\begin{equation}\label{Gauss}
\begin{aligned}
 |z|<1&;\\
 |z|=1&\text{ and } b_2-b_0-b_1>0;\\
 |z|=1&,\;z\ne1\text{ and } -1<b_2-b_0-b_1\le0.
\end{aligned}
\end{equation}
Here \[(b)_n=b\cdot(b+1)\cdot\ldots\cdot(b+n-1)=\Gamma(b+n)/\Gamma(b)\] is the Pochhammer symbol.
This function is defined by analytic continuation everywhere in $z,\;z<-1,$ see \cite{Andrews}. 

\subsection{Attracting-only case, $\gamma_0,\; \gamma_1>0.$}

 In this case, both parameters $\gamma_0$ and $\gamma_1$ are
regarded as positive revertive rates, and 
both patterns, $\phi_0$ and $\phi_1,$ defined by \eqref{def:phi}, 
converge as $t\to\infty,$
\[
\lim_{t\to\infty}\phi_0(t, x)=\frac{a_0}{\gamma_0}=:\rho_0,\qquad 
\lim_{t\to\infty}\phi_1(t, x)=\frac{a_1}{\gamma_1}=:\rho_1.
\]
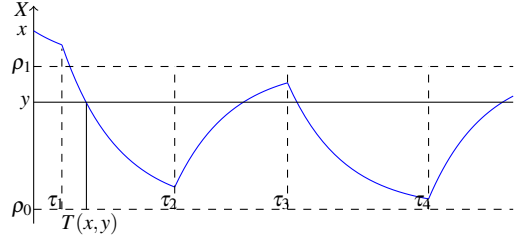
\begin{figure}[th]
 \begin{center}
\begin{tikzpicture}[x=0.75cm,y=0.95cm][domain=0:8]   
\draw[->] (-1,-1.2) -- (-1,1.85); 
\draw[-] (-0.065,-1) -- (-0.065,0.5);
 \draw[blue,domain=-1:-0.5] 
 	plot(\x,{1+0.5*exp(-\x-1)}); 
 \draw[blue,domain=-0.5:1.5] 
 plot(\x,{0.5*exp(-\x-1)-1+2*exp(-\x-0.5)}); 
  \draw[blue,domain=1.5:3.5] 
 plot(\x,{1+(-2+0.5*exp(-2.5)+2*exp(-2))*exp(-\x+1.5)}); 
 \draw[blue,domain=3.5:6] 
  plot(\x,{(1+(-2+0.5*exp(-2.5)+2*exp(-2))*exp(-4.5+2.5))*exp(-(\x-3.5))-1+exp(-\x+3.5)}); 
   \draw[blue,domain=6:7.5] 
 plot(\x,{((1+(-2+0.5*exp(-2.5)+2*exp(-2))*exp(-4.5+2.5))*exp(-(7-4.5))-1+exp(-7+4.5))*exp(-\x+6)+1-exp(-\x+6)}); 
 \draw [dashed] (-1,1)     to  (7.5,1); 
  \draw [dashed] (-1,-1)     to  (7.5,-1); 
  \draw (-1,0.5) to (7.5,0.5);
  \draw [dashed] (-0.5,-1)     to  (-0.5,1.25);   
  \draw [dashed] (1.5,-1)     to  (1.5,1); 
    \draw [dashed] (3.5,-1)     to  (3.5,1); 
      \draw [dashed] (6,-1)     to  (6,1); 
        \node at (-0.6,-0.9) {\scriptsize{$\tau_1$}};
                \node at (1.4,-0.9) {\scriptsize{$\tau_2$}};
                        \node at (3.4,-0.9) {\scriptsize{$\tau_3$}};
                                \node at (5.9,-0.9) {\scriptsize{$\tau_4$}};
    \node at (-1.2,1.8) {\scriptsize{$X$}};
        \node at (-1.2,1.5) {\scriptsize{$x$}};
 \node at (-1.2,1) {\scriptsize{$\rho_1$}};
  \node at (-1.2,-1) {\scriptsize{$\rho_0$}};
  \node at (-1.15, 0.5) {\scriptsize{$y$}};
  \node at (0,-1.2) {\scriptsize{$T(x, y)$}};

\end{tikzpicture}
 \caption{A sample path of $X=X(t),$ $\gamma_0, \gamma_1>0.
$ {}{}}
 \label{fig1}
 \end{center}
\end{figure} 
 %

The interval $[\rho_0,\;\rho_1]$
serves as an attractor for the paths of $X$: if the process starts at point $x$ outside this interval, $x\notin [\rho_0,\;\rho_1],$
 it falls into $[\rho_0,\;\rho_1]$ a.s. in a finite time. 
 Moreover, once caught, the process remains there forever,
 see \cite{MCAP2020}.  
  In this regard, we study the first passage through the threshold $y,\; \rho_0<y<\rho_1.$
 A sample path is shown in Fig. \ref{fig1}. 
 
The first passage time $T(x, y),$ \eqref{def:Txy}\textup, is finite a.s.  $\forall x,\;y,\;y\in(\rho_0,\;\rho_1).$
The distribution of $T(x, y)$ can be studied separately for $x<y$ and $x>y.$
In both these cases,  functions $\ell_0$ and $\ell_1$ corresponding to the Laplace transform
obey the boundary value problems for equations  \eqref{eq:hatQ++diff} on the half-lines 
 $x<y$ and $x>y,$ respectively.

First, let  $x<y,\;y\in(\rho_0,\;\rho_1).$ 

Since both levels, $\rho_0,\;\rho_1,$ attract, then
$t_0(x, y)=+\infty,\;t_1(x, y)=\frac{1}{\gamma_1}\log\frac{x-\rho_1}{y-\rho_1}<\infty.$ 
If the initial state is $1=\ep(0),$ then $T(x, y)\to0\;a.s.$ 
as $x\uparrow y.$
 The latter gives the boundary 
 condition
\begin{equation}
\label{eq:ell1=1}
\ell_1(q, x, y)|_{x\uparrow y}=1
\end{equation}
 to system \eqref{eq:hatQ++diff}. 
We solve this system on the half-line $x<y$ with boundary condition \eqref{eq:ell1=1}, writing the solution
in the form
\begin{equation}
\label{eq:ell0}
\begin{aligned}
\ell_0(q, x, y)=&\sum_{n=0}^\infty A_n(q, y)\xi_0(x)^n,\\
\ell_1(q, x, y)=&\sum_{n=0}^\infty B_n(q, y)\xi_0(x)^n.
\end{aligned}\end{equation}
Substituting functions $\ell_0$ and $\ell_1,$  defined by series \eqref{eq:ell0}, into \eqref{eq:hatQ++diff} 
and using the identities 
\begin{equation*}
(x-\rho_0)\frac{\rmd \xi_0(x)}{\rmd x}\equiv\xi_0(x),\qquad 
(\rho_1-\rho_0)\frac{\rmd \xi_0(x)}{\rmd x}\equiv1,
\end{equation*}
we obtain 
\begin{equation}
\label{eq:ABn}
\left\{
\begin{aligned}
 nA_n&=-\beta_0(q)A_n+\frac{\la_0}{\gamma_0}B_n,   \\
    nB_n-(n+1)B_{n+1}& =\frac{\la_1}{\gamma_1}A_n-\beta_1(q)B_n.
\end{aligned}
\right.
\end{equation}

After a simple algebra, see e.g. \cite{MCAP2020}, we find the solution of system \eqref{eq:ABn}: 
\begin{equation}\label{eq:AB0}
\left\{
\begin{aligned}
   A_n &=\frac{\la_0}{q+\la_0} \cdot\frac{(b_0)_n(b_1)_n}{(1+\beta_0)_nn!} B_0, \\
    B_n&= \frac{(b_0)_n(b_1)_n}{(\beta_0)_nn!} B_0,
\end{aligned}\qquad n\ge0.
\right.
\end{equation}
Here, recall,  $\beta_0=\beta_0(q)=(q+\la_0)/\gamma_0$ and $b_0,\;b_1$ are defined by \eqref{def:b0b1}.

Due to \eqref{Gauss}, series \eqref{eq:ell0} 
with coefficients $A_n$ and $B_n$ determined by \eqref{eq:AB0}
converge if $|\xi_0(x)|<1,$ that is, if
\begin{equation}
\label{eq:x0}
2\rho_0-\rho_1<x<\rho_1.
\end{equation}
Therefore, the Laplace transform of the first passage time $T(x, y)$ for $y\in(\rho_0,\;\rho_1),$ 
and $x$ satisfying \eqref{eq:x0}
is expressed in terms of the Gaussian hypergeometric series, \eqref{def:F},
\begin{equation}\label{eq:ell0ell1++}
\begin{aligned}
  \ell_0(q, x, y)  & =\frac{\la_0}{q+\la_0}F(b_0, b_1; 1+\beta_0; \xi_0(x))\cdot B_0,  \\
 \ell_1(q, x, y) &  =F(b_0, b_1; \beta_0; \xi_0(x))\cdot B_0.
\end{aligned}
\end{equation}

The indefinite parameter $B_0$ follows from the boundary condition \eqref{eq:ell1=1}. Thus, we finally obtain 
the explicit formulae for $\ell_0(q, x, y)$ and $\ell_1(q, x, y)$ in the case $\rho_0\le y<\rho_1$ 
and  $2\rho_0-\rho_1<x<y,$
\begin{equation}
\label{eq:ell0ell1++x<y}
\begin{aligned}
\ell_0(q, x, y)
=&\frac{\la_0}{q+\la_0}\cdot\frac{F(b_0, b_1;1+\beta_0(q); \xi_0(x))}{F(b_0, b_1; \beta_0(q); \xi_0(y))}, \\
\ell_1(q, x, y)
=&\frac{F(b_0, b_1; \beta_0(q); \xi_0(x))}{F(b_0, b_1; \beta_0(q); \xi_0(y))}.
\end{aligned}
\end{equation}

In the case  $\rho_0<y\le\rho_1,\;x>y$ formulae for $\ell_0(q, x, y)$ and $\ell_1(q, x, y)$
can be obtained by symmetry in the form of series with $\xi_1(x).$
We have $\ell_0(q, x, y)|_{x\downarrow y}=1$  and 
\begin{equation}
\label{eq:ell0ell1++x>y}
\begin{aligned}
\ell_0(q, x, y)&
=\frac{F(b_0, b_1; \beta_1(q); \xi_1(x))}{F(b_0, b_1; \beta_1(q); \xi_1(y))},\\
\ell_1(q, x, y)&
=\frac{\la_1}{q+\la_1}\cdot\frac{F(b_0, b_1;1+\beta_1(q); \xi_1(x))}{F(b_0, b_1; \beta_1(q); \xi_1(y))},
\end{aligned}
\end{equation}
\[\rho_0<y<x<2\rho_1-\rho_0,\qquad y\le\rho_1.\]

Note that for $y\in(\rho_0,\;\rho_1)$ we have $0<\xi_0(y),\xi_1(y)<1,$ therefore the denominators in
\eqref{eq:ell0ell1++x<y} and \eqref{eq:ell0ell1++x>y} are also determined by Gaussian hypergeometric series \eqref{def:F}.

In the case of non-strict attraction, the Laplace transforms $\ell_0$ and $\ell_1$
can be obtained similarly to \eqref{eq:ell0ell1++x<y} and \eqref{eq:ell0ell1++x>y}. 
For example, let $\gamma_1=0, a_1=+1,\;\gamma_0>0$ and $x<y.$ We have:
\[
\begin{aligned}
    \ell_0(q, x, y)=&\frac{\la_0/\gamma_0}{\beta_0(q)}
    \frac{\Phi(\delta; 1+\beta_0(q); (x-\rho_0)(q+\la_1))}{\Phi(\delta; \beta_0(q); (y-\rho_0)(q+\la_1)))},   \\
  \ell_1(q, x, y)=  &    \frac{\Phi(\delta; \beta_0(q); (x-\rho_0)(q+\la_1))}{\Phi(\delta; \beta_0(q); (y-\rho_0)(q+\la_1)))},
\end{aligned}
\]
where $\delta=\dfrac{(q+\la_0)(q+\la_1)-\la_0\la_1}{\gamma_0(q+\la_1)}$ 
and $\Phi(\cdot;\cdot;\cdot)$ is the confluent hypergeometric function.

\subsection{Attraction-repulsion:
$\gamma_0$ and $\gamma_1$ have opposite signs, ``raznotyk''}
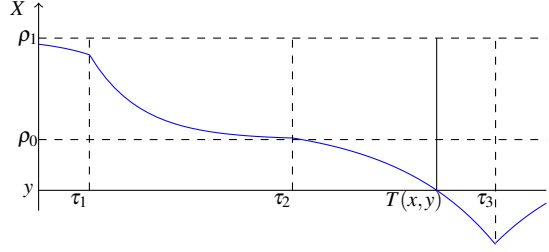
\begin{figure}[th]
 \begin{center}
\begin{tikzpicture}[x=1.35cm,y=1.35cm][domain=0:8]   
\draw[->] (-1,-0.2) -- (-1,1.85); 
\draw[-] (2.92,0) -- (2.92,1.5);
 \draw [dashed] (-1,1.5)     to  (4,1.5); 
  \draw [dashed] (-1,0.5)     to  (4,0.5); 
    \draw [] (-1,0)     to  (4,0); 
 \draw[blue,domain=-1:-0.5] 
 	plot(\x,{1.5-2*exp(2*\x-1.5)}); 
 \draw[blue,domain=-0.5:1.5] 
 plot(\x,{0.5+2*exp(-2*\x-1.87)}); 
  \draw[blue,domain=1.5:3.5] 
plot(\x,{0.68-2*exp(\x-4)});  
 \draw[blue,domain=3.5:4] 
  plot(\x,{0.48-exp(-\x+3.5)}); 

  \draw [dashed] (-0.5,0)     to  (-0.5,1.5);   
  \draw [dashed] (1.5,0)     to  (1.5,1.5); 
    \draw [dashed] (3.5,-0.48)     to  (3.5,1.5); 
        \node at (-0.6,-0.1) {\scriptsize{$\tau_1$}};
                \node at (1.4,-0.1) {\scriptsize{$\tau_2$}};
                        \node at (3.4,-0.1) {\scriptsize{$\tau_3$}};
    \node at (-1.2,1.8) {\scriptsize{$X$}};
        \node at (-1.1,1.5) {\scriptsize{$\rho_1$}};
        \node at (-1.1,0.5) {\scriptsize{$\rho_0$}};
                \node at (-1.1,0) {\scriptsize{$y$}};
  \node at (2.7,-0.1) {\scriptsize{$T(x,y)$}};
%
\end{tikzpicture}
 \caption{The sample path of $X=X(t),$ $\gamma_0>0>\gamma_1. 
$ {}{}}
 \label{fig2}
 \end{center}
\end{figure}

To be specific, assume that $\gamma_0>0>\gamma_1,$ that is,
the pattern $\phi_1$ is repelled from the threshold $\rho_1,$ while $\phi_0$ is attracted to $\rho_0$.
In this case, 
process $X=X(t),\;t\ge0,$ a.s. falls under $\rho_0,$ into the half-line $\{z~|~z\le\rho_0\},$
in a finite time and, once falling, remains there forever, see Fig. \ref{fig2}.
Let the threshold $y$ belong to the attractor, $y<\rho_0.$ Similarly to the case of two attractive levels, 
we obtain the boundary conditions in dependence of the starting point $x.$ 

%

If the process begins with $0$-state from below of threshold $y,\;x<y,$  then $T(x, y)\to0,\;a.s.$ as $x\uparrow y.$
For $x<y,$ we must consider the system \eqref{eq:hatQ++diff}  with the boundary condition 
\begin{equation}
\label{eq:ell0=1}
\ell_0(q, x, y)|_{x\uparrow y}=1,
\end{equation}
which gives the solution of the form \eqref{eq:ell0ell1++}, $2\rho_0-\rho_1<x<y<\rho_0.$
By the boundary condition \eqref{eq:ell0=1} we obtain 
\[
B_0=\left(\frac{\la_0}{q+\la_0}F(b_0, b_1; 1+\beta_0; \xi_0(x))\right)^{-1}
\]
and the explicit solution is given by
\begin{equation*}
\begin{aligned}
\ell_0(q, x, y)=&\frac{F(b_0, b_1; 1+\beta_0; \xi_0(x))}{F(b_0, b_1; 1+\beta_0; \xi_0(y))},\\
\ell_1(q, x, y)=&\frac{q+\la_0}{\la_0}\cdot\frac{F(b_0, b_1; \beta_0; \xi_0(x))}{F(b_0, b_1; 1+\beta_0; \xi_0(y))},
\end{aligned}\end{equation*}
\[2\rho_0-\rho_1<x<y<\rho_0.\]

Acting in the similar way, we find that in the case $2\rho_0-\rho_1<y<x<\rho_1,\;y<\rho_0$  we obtain,
\begin{equation*}
\begin{aligned}
\ell_0(q, x, y)=&\frac{q+\la_1}{\la_1}\cdot\frac{F(b_0, b_1; \beta_1; \xi_1(x))}{F(b_0, b_1; 1+\beta_1; \xi_1(y))},\\
\ell_1(q, x, y)=&\frac{F(b_0, b_1; 1+\beta_1; \xi_1(x))}{F(b_0, b_1; 1+\beta_1; \xi_1(y))}.
\end{aligned}\end{equation*}

In the case of $\gamma_0<0<\gamma_1,$ 
the formulae for the distribution of $T(x, y),\;y>\rho_1,$
are symmetric.
\section{Invariant measures}\label{sec5}

As well as the distribution of the first passage time, 
the form of the invariant measure differs 
in the case of positive values of $\gamma_0,\;\gamma_1$ 
and in the case of opposite signs.
\subsection{Attracting-only case, $\gamma_0,\; \gamma_1>0$}\label{sec41}

Since the paths of $X$ remain inside the interval $(\rho_0,\;\rho_1)$ after an almost surely finite transition time, 
the invariant measure $\vec\mu$ is supported on 
$[\rho_0,\;\rho_1],\;\rho_0=a_0/\gamma_0,$ $\rho_1=a_1/\gamma_1.$

As it was shown in Section \ref{sec2},
the invariant probability density function 
$\vec\pi=\vec\pi(x)=(\pi_0(x),\;\pi_1(x))$ obeys the system \eqref{eq:L*system}
of the ordinary differential equations, $\rho_0<x<\rho_1$.
In the attracting-only case,
this system should be equipped with the boundary conditions, see \eqref{eq:L*bc},
\begin{equation}
\label{eq:boundary++}
\begin{aligned}
\pi_0(x)|_{x=\rho_1-}=0,&\quad  (x-\rho_0)\pi_0(x)|_{x=\rho_0+}=0,\\
(\rho_1-x)\pi_1(x)|_{x=\rho_1-}=0,&\quad  \pi_1(x)|_{x=\rho_0+}=0.
\end{aligned}
\end{equation}

By substituting 
\[\begin{aligned}
\pi_0(x)=&C_0(x-\rho_0)^{k_{00}}(\rho_1-x)^{k_{01}}\1_{\{\rho_0<x<\rho_1\}},\\
\pi_1(x)=&C_1(x-\rho_0)^{k_{10}}(\rho_1-x)^{k_{11}}\1_{\{\rho_0<x<\rho_1\}},
\end{aligned}\]
into equations \eqref{eq:L*system} and taking into account the boundary conditions \eqref{eq:boundary++},
we obtain
\begin{equation*}
\begin{aligned}
  \pi_0(x)=  & c\gamma_0^{-1}(x-\rho_0)^{-1+\alpha_0}(\rho_1-x)^{\alpha_1}  \1_{\{\rho_0<x<\rho_1\}},  \\
   \pi_1(x)=  &c\gamma_1^{-1} (x-\rho_0)^{\alpha_0}(\rho_1-x)^{-1+\alpha_1}  \1_{\{\rho_0<x<\rho_1\}},  
\end{aligned}
\end{equation*}
where $\alpha_0=\la_0/\gamma_0, \alpha_1=\la_1/\gamma_1,\;\alpha_0, \alpha_1>0.$ 
The normalising constant $c$ can be found from 
the equality \[\int_{\rho_0}^{\rho_1}[\pi_0(x)+\pi_1(x)]\rmd x=1.\]
Due to \cite[3.196]{GR} we obtain
\begin{equation}
\label{eq:c-inv}
\begin{aligned}
c^{-1}
=&(\rho_1-\rho_0)^{\alpha_0+\alpha_1}
\Big[\gamma_0^{-1}B(\alpha_0,\;\alpha_1+1))\\
+&\gamma_1^{-1}B(\alpha_0+1,\;\alpha_1))\Big],
\end{aligned}\end{equation}
where $B(\cdot,\;\cdot)$ is Euler's beta-function.

\begin{ex}
Let $\gamma_0=\gamma_1=\gamma>0.$ 
 In this case\textup, by \eqref{eq:c-inv}
$
c^{-1}=(\rho_1-\rho_0)^2\gamma^{-1}
$
and 
\[\begin{aligned}
\pi_0(x)=&(\rho_1-\rho_0)^{-2}(\rho_1-x),\\ \pi_1(x)=&(\rho_1-\rho_0)^{-2}(x-\rho_0),
\end{aligned}\]
\[
\qquad \rho_0<x<\rho_1.
\]
It's curious that when the process $X$ begins with states 0 or 1 with equal probability\textup,
 the invariant distribution is uniform on $[\rho_0,\;\rho_1]$.
\end{ex}

\subsection{Attraction-repulsion: 
$\gamma_0$ and $\gamma_1$ have opposite signs, ``raznotyk''}
First, let $\gamma_0>0>\gamma_1.$ In this case, after an almost surely finite transition period,
process $X$ falls into the half-line $\{x<\rho_0\},$ see Fig. \ref{fig2}.
The invariant distributions are defined by the probability density functions $\vec\pi=(\pi_0(x),\;\pi_1(x)),$
satisfying system \eqref{eq:L*system}, $x<\rho_0,$ with the boundary conditions
\[
(\rho_0-x)\pi_0(x)|_{x=\rho_0}=0,\qquad
(\rho_1-x)\pi_1(x)|_{x=\rho_0}=0.
\]
Similarly the attraction-only case, Section \ref{sec41},
 the solution is given by 
\begin{equation}
\label{eq:pi0pi1+-}
\begin{aligned}
   \pi_0(x) &=C_0(\rho_0-x)^{-1+\alpha_0}(\rho_1-x)^{\alpha_1}\1_{\{x<\rho_0\}},   \\
    \pi_1(x)&  =C_1(\rho_0-x)^{\alpha_0}(\rho_1-x)^{-1+\alpha_1}\1_{\{x<\rho_0\}},
\end{aligned}
\end{equation}
where $C_0=c_0\gamma_0^{-1}, C_1=-c_0\gamma_1^{-1},\;C_0, C_1>0,$ and 
$\alpha_0=\la_0/\gamma_0, \alpha_1=\la_1/\gamma_1,\;\alpha_0>0>\alpha_1.$

In the symmetric case, $\gamma_0<0<\gamma_1,$ process $X$ is captured by the upper half-line, $K=\{x>\rho_1\}.$
The invariant probability density function is determined by
\begin{equation}
\label{eq:pi0pi1-+}
\begin{aligned}
   \pi_0(x) &=C_0(x-\rho_0)^{-1+\alpha_0}(x-\rho_1)^{\alpha_1}\1_{\{x>\rho_1\}},   \\
    \pi_1(x)&  =C_1(x-\rho_0)^{\alpha_0}(x-\rho_1)^{-1+\alpha_1}\1_{\{x>\rho_1\}},
\end{aligned}
\end{equation}
where $C_0=-c_1\gamma_0^{-1},\;C_1=c_1\gamma_1^{-1}.$

In both cases, we assume that 
\begin{equation}
\label{eq:alpha<0}
\alpha_0+\alpha_1<0.
\end{equation}
The normalising constants $c_0$ and $c_1$ are determined by the condition
\begin{equation}
\label{eq:1}
\int_K\left[\pi_0(x)+\pi_1(x)\right]\rmd x=1.
\end{equation}
The integral in \eqref{eq:1} converges if \eqref{eq:alpha<0} holds.
In the case $\gamma_0>0>\gamma_1,$ the normalising constant 
$c_0$ for \eqref{eq:pi0pi1+-} is determined by
\begin{equation*}
\begin{aligned}
c_0^{-1}
=(\rho_1-\rho_0)^{\alpha_0+\alpha_1}&
\big(\gamma_0^{-1}B(-\alpha_0-\alpha_1, \alpha_0)\\
-&\gamma_1^{-1}B(-\alpha_0-\alpha_1, 1+\alpha_0)\big), 
\end{aligned}\end{equation*}
\[\alpha_0+\alpha_1<0,\quad \alpha_0>0.\]
In the symmetric case $\gamma_0<0<\gamma_1,$ 
the normalising constant $c_1$ for \eqref{eq:pi0pi1-+} is determined by
\begin{equation*}
\begin{aligned}
c_1^{-1}
=(\rho_1-\rho_0)^{\alpha_0+\alpha_1}&
\big(-\gamma_0^{-1}B(-\alpha_0-\alpha_1, 1+\alpha_1)\\
+&\gamma_1^{-1}B(-\alpha_0-\alpha_1, \alpha_1)\big).
\end{aligned}\end{equation*}
\[\alpha_0+\alpha_1<0,\quad \alpha_1>0.\]

If, on the contrary, \eqref{eq:alpha<0} is not met, the invariant probability distribution does not exist.

Repulsion-only case, $\gamma_0, \gamma_1<0,$ corresponds to 
a subordinator, that is, all paths of $X(t),\;t\ge0$ (or $-X(t),$ $t\ge0$) are strictly monotonically increasing.
Therefore, there are no invariant probability distributions.


\subsection{Non-strict attraction}
If one of attraction rates is zero, say $\gamma_1=0,$ and $\gamma_0>0,$ then the pattern
$\phi_0$ defined by \eqref{def:phi} is attractive, and $\phi_1(t, x)\equiv x+a_1t.$ 

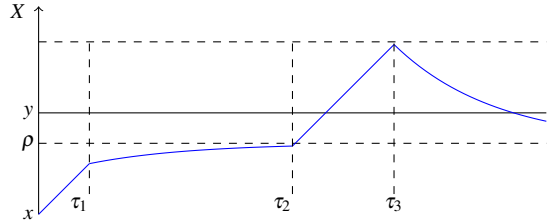
\begin{figure}[th]
 \begin{center}
\begin{tikzpicture}[x=1.35cm,y=1.35cm][domain=0:8]   
\draw[->] (-1,-0.2) -- (-1,1.85); 
 \draw [dashed] (-1,1.5)     to  (4,1.5); 
  \draw [dashed] (-1,0.5)     to  (4,0.5); 
    \draw [] (-1,0.8)     to  (4,0.8); 
 \draw[blue,domain=-1:-0.5] 
 	plot(\x,{0.8+\x});      
 \draw[blue,domain=-0.5:1.5] 
 plot(\x,{0.5-2*exp(-\x-2.80)}); 
  \draw[blue,domain=1.5:2.5] 
plot(\x,{0.473+\x-1.5});  
 \draw[blue,domain=2.5:4] 
  plot(\x,{0.5+2*exp(-\x+1.78)}); 

  \draw [dashed] (-0.5,0)     to  (-0.5,1.5);   
  \draw [dashed] (1.5,0)     to  (1.5,1.5); 
    \draw [dashed] (2.5,0)     to  (2.5,1.5); 
        \node at (-0.6,-0.1) {\scriptsize{$\tau_1$}};
                \node at (1.4,-0.1) {\scriptsize{$\tau_2$}};
    \node at (-1.2,1.8) {\scriptsize{$X$}};
        \node at (-1.1,0.5) {\scriptsize{$\rho$}};
                \node at (-1.1,0.8) {\scriptsize{$y$}};         
                \node at (-1.1,-0.2) {\scriptsize{$x$}};
  \node at (2.5,-0.1) {\scriptsize{$\tau_3$}};
%
\end{tikzpicture}
 \caption{The sample path of $X=X(t),$ $\gamma_0>0=\gamma_1,$ $a_1=1. 
$ {}{}}
 \label{fig3}
 \end{center}
\end{figure} 

Without loss of generality, we will put $a_1=1$ (if, originally, $a_1$ is positive), 
 $a_1=-1$  (if, originally, $a_1$ is negative). 
 Regardless of the initial point $x,$ the trajectories of $X$ 
(possibly, after a finite transition period) remain above (or below) the threshold $\rho=\rho_0=a_0/\gamma_0$ 
in the case $a_1=+1$ (or, respectively, $a_1=-1$), see Fig.\ref{fig3}.

The invariant distribution densities $\pi_0(x),\;\pi_1(x)$ follow system \eqref{eq:L*system} 
on $x<\rho,$ if $a_1=-1$ (on $x>\rho$, if $a_1=+1$): 
\[\left\{
\begin{aligned}
    (x-\rho)\frac{\rmd \pi_0(x)}{\rmd x}&=\left(\frac{\la_0}{\gamma_0}-1\right)\pi_0(x)-\frac{\la_1}{\gamma_0}\pi_1(x),   \\
   -a_1\frac{\rmd \pi_1(x)}{\rmd x} &=-\la_0\pi_0(x)+\la_1\pi_1(x).  
\end{aligned}\right.
\]
In both cases, the boundary condition \eqref{eq:L*bc} turns into
\begin{equation}\label{eq:bc00}
(\rho-x)\pi_0(x)|_{x=\rho}=0,\qquad\pi_1(x)|_{x=\rho}=0.
\end{equation}
The solutions of these two boundary value problems are given by
\begin{equation}
\label{eq:pi0}
\begin{aligned}
    \pi_0^\pm(x)&=C|x-\rho|^{\alpha-1}\rme^{-\la_1|x-\rho|}\theta^\pm(x),  \\
    \pi_1^\pm(x)&=\gamma_0C |x-\rho|^{\alpha}\rme^{-\la_1|x-\rho|}\theta^\pm(x),
\end{aligned}
\end{equation}
where $\alpha=\la_0/\gamma_0,\; \rho=a_0/\gamma_0;$  
$C$ is a normalising constant.
Functions $\theta^-(x)=\1_{\{x<\rho\}}$ and $\theta^+(x)=\1_{\{x>\rho\}}$ 
specify the measure support below and, respectively, above the level $\rho.$
The explicit value of normalising constant $C$ follows from the equalities
\[\begin{aligned}
1=&\int_{\rho}^\infty[\pi_0^+(x)+\pi_1^+(x)]\rmd x\\
=&C\la_1^{-\alpha}\Gamma(\alpha)+\gamma_0C\la_1^{-(\alpha+1)}\Gamma(\alpha+1)\\
=&C\la_1^{-\alpha}\Gamma(\alpha)\left[1+\la_0/\la_1\right],
\end{aligned}\]
which gives 
\[
C=\frac{\la_1^{1+\alpha}}{(\la_0+\la_1)\Gamma(\alpha)}.
\]
The proof for the case $a_1=-1$ is similar.

\begin{rem}
When $\gamma_1=0$ and $a_1=0,$ system \eqref{eq:L*system} turn into
\[
\begin{aligned}
 (x-\rho)\pi_0'(x)=&-\pi_0(x),  \\
0=&-\la_0\pi_0(x)+\la_1\pi_1(x),
\end{aligned}
\]
 with the boundary conditions \eqref{eq:bc00}.
This means $\pi_0=\pi_1=0,$ i.e. in this case there is no invariant probability measure.
\end{rem}
\section{Scaling}\label{sec6}

First, we speed up the underlying Markov chain. Let $\ep$ be driven by 
alternating switching intensities $\la_0,\;\la_1$ which are 
 high but comparable, i.e. 
\begin{equation}
\label{eq:nu}
\la_0, \la_1\to\infty\qquad\text{and}\qquad \frac{\la_0}{\la_1}\to\nu,\qquad \nu>0.
\end{equation}
With this scaling, the invariant distribution of the underlying process $\ep,\;\pi_*(\la_0,\;\la_1),$ 
\eqref{eq:pi*}, becomes
\[
\pi_*(\la_0,\;\la_1)\to\pi_*^*=\left((1+\nu)^{-1},\;\nu(1+\nu)^{-1}\right).
\]
Passing to limit in the moment generating function, see \cite[formula (3.3) ]{STAPRO82},
it is easy to see that if the velocities $c_0,\;c_1$  remain frozen, 
then upon scaling \eqref{eq:nu}
the telegraph process $\mathbb T(t),\;t>0,$ \eqref{def:telegraph},
converges in probability to $c_*t$ uniformly in $ t\in[0, T],$
\begin{equation*}
\mathbb T(t)=\int_0^tc_{\ep(s)}\rmd s\to c_*t,
\end{equation*}
where $c_*=\vec\pi_*^*\cdot\vec c.$ 

This result help to understand the limiting behaviour of the Markov-modulated Ornstein-Uhlenbeck process.
Let the switching intensities tend to infinity, \eqref{eq:nu}, and 
parameters $\vec a=(a_0, a_1),\;\vec b=(b_0, b_1),\;\vec\gamma=(\gamma_0, \gamma_1)$
of process $M,$ \eqref{eq:solKacLangevind}, remain constant.
With this scaling, process $M$  weakly converges to process $M_*,$
which is an ordinary (unmodulated) Ornstein-Uhlenbeck process
 defined by the Langevin equation, \eqref{eq:KacLangevind}, 
with constant \emph{deterministic parameters} 
$a_\infty=\vec\pi_*^*\vec a^{\mathrm T},\;b_\infty=\vec\pi_*^*\vec b^{\mathrm T},$ 
and $\gamma_\infty=\vec\pi_*^*\vec \gamma^{\mathrm T},$ 
\[
\rmd M_*(t)= \left(a_\infty-\gamma_\infty M_*(t)\right)\rmd t+b_\infty\rmd W(t),\qquad t>0.
\]
 See \cite[Corollary 5.1]{huang2014}.
 
 Further, let's see what happens under different versions of  Kac scaling, see \cite{Kac}.

It is known that with the classical Kac scaling, i.e. if
\begin{equation}
\label{eq:kac}
\begin{aligned}
c_0,\;-c_1\to+\infty,&\qquad \la_0,\;\la_1\to+\infty,\\
\text{ and }&~~\\
c_0^2/\la_0,\qquad & c_1^2/\la_1\to\sigma^2,
\end{aligned}\end{equation}
 the telegraph process 
$\mathbb T,$ \eqref{def:telegraph}, converges in distribution 
on $C\left([0, T];\RR\right)$ (equipped with the sup-norm)
  to the Wiener process $\sigma W(t),$ see the proof e.g. in
 \cite{MPRF}. 
 
 We apply this idea to a substantially \emph{asymmetric} telegraph process. 
Let the scaling condition similar to \eqref{eq:kac}
be satisfied separately for the two states, i. e. let \eqref{eq:nu} holds, and
$c_0\to+\infty$, $c_1\to-\infty,$ so that
\begin{equation}
\label{Kac-inhom}
\frac{c_0}{\sqrt{\lambda_0}}\to\sigma_0,\qquad
\frac{c_1}{\sqrt{\lambda_1}}\to-\sigma_1,
\end{equation}
where 
$\sigma_0, \sigma_1>0$.
From \eqref{eq:nu}-\eqref{Kac-inhom} it follows that the velocities are also comparable:
\[
\frac{c_0}{c_1}=\frac{c_0/\sqrt{\lambda_0}}{c_1/\sqrt{\lambda_1}}
\sqrt{\lambda_0/\lambda_1}\to -\nu\sigma_0/\sigma_1.
\]
We also assume that 
\begin{equation}\label{eq:delta}
\frac{\lambda_1c_0+\lambda_0c_1}{\lambda_1+\lambda_0}\to\delta.
\end{equation}
The latter limit relation is equivalent to 
$
c_1\dfrac{c_0/c_1+\lambda_0/\lambda_1}{1+\lambda_0/\lambda_1}\to\delta.
$ 
Therefore, condition \eqref{eq:delta} assumes that
$c_0/c_1+\lambda_0/\lambda_1\to0.$ 
So, condition \eqref{eq:delta} reads as rate of  ``similarity" between 
$\lambda_0/\lambda_1$ and $c_0/c_1$ at infinity. More precisely, 
\[
c_1(c_0/c_1+\lambda_0/\lambda_1)\to\delta(1+\nu^2).
\]

Under the scaling conditions \eqref{eq:nu}-\eqref{eq:delta} stated above,
the telegraph process $\mathbb T(t)$ weakly converges to the Wiener process with drift, 
see \cite{JAP51},
\begin{equation}
\label{eq:convergenceKac}
\mathbb T(t)=\int_0^tc_{\ep(s)}\rmd s\stackrel{\mathcal D}{\to} \sigma W(t)+\delta t,\qquad t>0.
\end{equation}
Here 
\begin{equation}\label{def:sigma}
\sigma=\frac{\sigma_0\sigma_1}{\sqrt{(\sigma_0^2+\sigma_1^2)/2}}.
\end{equation}

 
Assuming \eqref{eq:nu} to be hold, we apply scaling conditions similar to
\eqref{Kac-inhom}-\eqref{eq:delta} in several versions. 

 \begin{list}{\textbf{({\alph{lis}})}}
{\usecounter{lis}}
\item\label{a}
Let conditions \eqref{Kac-inhom}-\eqref{eq:delta} be satisfied for $a_0,\; a_1$. That is, let 
$a_0\to-\infty, \; a_1\to+\infty,$ so that 
$a_0/\sqrt{\la_0}\to-\sigma_{0,\;a},\; a_1/\sqrt{\la_1}\to\sigma_{1,\;a}.$

Let the additional drift be caused by condition of the form \eqref{eq:delta},
\begin{equation*}
\frac{\la_1a_0+\la_0a_1}{\la_0+\la_1}\to\delta_a.
\end{equation*} 

Leaving the remaining parameters constant, we see that due to \eqref{eq:solKacLangevind},
 process $M(t),$ \eqref{eq:KacLangevind}, weakly converges to $M_a(t),$
 \[
 M_a(t)=\rme^{-\gamma_\infty t}\left(
x+\int_0^t\rme^{\gamma_\infty s} \rmd\left(\sigma_a\widetilde W(s)+\delta_as\right)
+\int_0^t\rme^{\gamma_\infty s}b_{\ep(s)}\rmd W(s)
 \right),\qquad t>0,
 \]
which is the solution of the equation
\begin{equation}
\label{eq:M-A}
\rmd M_a(t)=\left(\delta_a-\gamma_\infty M_a(t)\right)\rmd t+\sigma_a\rmd\widetilde W(t)+b_\infty\rmd W(t),
\end{equation}
where $\widetilde W$ is the Wiener process independent of $W,$
and $\sigma_a$ is defined by formula \eqref{def:sigma} 
with $\sigma_0=\sigma_{0, a}$ and $\sigma_1=\sigma_{1, a}.$

\item
Assume $\vec a$ and $\vec b$ to be constant, but
 conditions \eqref{Kac-inhom}-\eqref{eq:delta} be satisfied for $\gamma_0,\; \gamma_1:$
$
\gamma_0\to-\infty, \; \gamma_1\to+\infty,
$
with the consistency conditions
\begin{equation*}
\frac{\gamma_0}{\sqrt{\la_0}}\to-\sigma_{0,\;\gamma},\qquad \frac{\gamma_1}{\sqrt{\la_1}}\to\sigma_{1,\;\gamma}
\end{equation*}
and additional drift caused by
\begin{equation*}
\frac{\la_1\gamma_0+\la_0\gamma_1}{\la_0+\la_1}\to\delta_\gamma.
\end{equation*} 

With this scaling, similarly to \eqref{eq:M-A}, 
we obtain the limiting process $M_b$ satisfying the equation
\begin{equation*}
\rmd M_b(t)=\left(a_\infty-\delta_\gamma M_b(t)\right)\rmd t
-\sigma_\gamma M_b(t)\rmd \widetilde W(t)
+b_\infty\rmd W(t),\quad t>0,
\end{equation*}
where $\widetilde W$ is the Wiener process independent of $W$ and $\sigma_\gamma$ is defined 
by formula \eqref{def:sigma} 
with $\sigma_0=\sigma_{0, \gamma}$ and $\sigma_1=\sigma_{1, \gamma}.$

\item
Let both conditions for $\vec a$, \textbf{(a)}, and conditions for $\vec\gamma$, \textbf{(b)}
 hold simultaneously,
and $\vec b$ be constant.

Then $M(t)$ converges to $M_{c}(t)$ determined by the equation, $t>0,$
\begin{equation*}
\rmd M_{c}(t)=\left(\delta_a-\delta_\gamma M_c(t)\right)\rmd t
+\left(\sigma_a-\sigma_\gamma M_c(t)\right)\rmd\widetilde W(t)
+b_\infty\rmd W(t).
\end{equation*}
\end{list}

\section*{Acknowledgements} 
This research was supported by the Russian Foundation for Basic Research (RFBR) 
and Chelyabinsk Region, project number 20-41-740020.


\end{document}